\documentclass[reqno,11pt]{amsart}
\usepackage{graphicx}
\usepackage{amsaddr}
\usepackage{hyperref}
\usepackage{etoolbox}
\apptocmd{\thebibliography}{\raggedright}{}{}
\usepackage{amsthm}
\usepackage{amssymb}
\usepackage{amsmath}
\usepackage{bbm}
\usepackage{url}
\usepackage[parfill]{parskip}
\usepackage{caption}
\usepackage{subcaption}
\usepackage{enumitem}

\newcommand{\rd}{\mathbb{R}^d}
\newcommand{\expect}{\mathbb{E}}

\newtheorem{innercustomthm}{Lemma}
\newenvironment{customthm}[1]
  {\renewcommand\theinnercustomthm{#1}\innercustomthm}
  {\endinnercustomthm}

\title{Explicit bounds for a Gaussian decomposition Lemma of Sellke}
\author{Tobias Schmidt }
\address[A1]{TU Darmstadt}
\email[A1]{tobias.schmidt@tu-darmstadt.de}
\begin{document}
\maketitle
\begin{abstract}
    In \cite{Se22} a decomposition of Gaussian measures on finite- \linebreak dimensional spaces was introduced, which turned out to be a central technical tool to improve currently known bounds on a long standing conjecture in statistical mechanics called the Polaron problem. This note slightly generalizes this decomposition and provides numerical values for all occurring constants.
\end{abstract}

The Gaussian correlation inequality (GCI) was proven recently by Royen \cite{Ro14} and is a central ingredient in \cite{Se22}, where substantial progress was made on a long-standing conjecture called the \textit{Polaron problem}. In order to apply GCI in a sensible manner it is shown in \cite[Lemma 3.1]{Se22} that any centred Gaussian measure $\mu$ on a finite dimensional, real normed space $\Omega$ can be written as a mixture of two measures $\nu_{\text{good}}$ and $\nu_{\text{bad}}$. $\nu_{\text{good}}$ is supported on an (essentially) compact set and is dominated by $\mu$ in the sense defined below. Crucially, $\nu_{\text{bad}}$ is dominated by a dilated version of $\mu$, allowing for easy error estimates on the part of $\mu$ where large deviations are possible. 
In this note we present a slightly more general, quantitative version of this mixture decomposition.

For two probability measures $\nu$ and $\mu$ we write $\nu \preceq \mu$ whenever $\frac{\mathrm{d}\nu}{\mathrm{d}\mu}$ is the finite product of non-negative, symmetric quasi-concave functions or the uniformly bounded limit of such functions. The relation $\preceq$ is interesting because whenever $\mu$ is centred Gaussian it is a consequence of GCI that
$$\nu(f) \leq \mu(f)$$
for $f$ quasi-convex and bounded. For more information on this form of Gaussian domination see e.g. \cite{Se22}.

Define for each Borel set $B \in \mathcal{B}(\Omega)$ and $n \geq 1$
$$\mu^{\times n}(B) := \mu (B/n).$$
It is easily verified that whenever $\mu$ is standard normal on $\rd$, then
$$\frac{\mathrm{d}\mu}{\mathrm{d}\mu^{\times n}} \propto \exp \left( - \frac{n^2 -1}{ 2 n^2} \Vert x \Vert^2 \right).$$ 
In the following, $K^c$ denotes the complement of $K$.
\begin{customthm}{1}
    \label{31_thm}
    Let $\mu$ be a centered Gaussian measure on $(\Omega$,$\mathcal{B}(\Omega))$. Furthermore, let $K \in \mathcal{B}(\Omega)$ be a symmetric, closed convex set with  $\mu(K^c) \leq \delta < 0.5$  and
$$\pi \left( 2\log(\delta^{-1})+8 \right) \leq \delta^{-1}.$$ 
Then, there exist two probability measures  $\nu_{good}, \nu_{bad} $ and $\delta'>0$ such that
    $$\mu = (1-\delta')\nu_{good} + \delta' \nu_{bad}.$$
Moreover, for $n > 1$, $C:=  32 \frac{n^2}{n^2 -1}$ it holds that
\begin{enumerate}[label=\textbf{R.\arabic*}]
\item $\delta' \leq 2\delta$\label{r1}, 
\item supp($\nu_{good}$) $ \subseteq C K$\label{r2},
\item $\nu_{good} \preceq \mu$\label{r3},
\item $\nu_{bad} \preceq \mu^{ \times n}$\label{r4}.
\end{enumerate}

\end{customthm}

We briefly compare the results from Lemma \ref{31_thm} with the original statement \cite[Lemma 3.1]{Se22}. The main difference is that the original statement shows the existence of a constant $C$, whereas we provide exact values for $\delta$ small enough.
Moreover, we allow arbitrary dilations of $\mu$ to dominate $\nu_\text{bad}$. Compared to the original result we only show $\delta' \leq 2 \delta$, but the reader will easily verify that the additional factor of $2$ can be omitted by increasing the value of $C$ accordingly. The same arguments can be used to trade the requirement 
$$\pi \left( 2\log(\delta^{-1})+8 \right) \leq \delta^{-1}$$
for just $\delta < 0.5$ (this condition suffices to always find a centred ball $B(\epsilon)$ of some radius $\epsilon > 0$ contained in $K$), but paying the price of increasing $C$. Indeed, scaling $K$ by an additional factor $\Tilde{c}$ until $\mu ((\Tilde{c}K)^c)$ is small enough for the condition in the last display to hold suffices.

We follow the original proof provided in \cite{Se22} for the most part and only change the arguments whenever it is more convenient with respect to $C$. For example, we do not provide a new argument to show \ref{r4}, but to show \ref{r1} we refrain from using the isoperimetric inequality in order to keep track of $C$. 

\begin{proof}[Proof of Lemma \ref{31_thm}]
    W.l.o.g. we may assume that $\Omega = \rd$ with $d=\dim(\Omega)$ and that $\mu$ is standard normal on $\rd$. We first claim that
    \begin{equation}
        \label{ball_contained}
        B(\sqrt{\log(\delta^{-1})})\subseteq K.
    \end{equation}
    To show this, recall
    that the standard normal on $\rd$ is rotation-invariant. Assume that there is a point $x \in K^c$ that is $c := \sqrt{\log(\delta^{-1})}$-close to the origin (i.e., the biggest radius ball we could fit into $K$ would be exactly that). Let $G$ be the hyperplane separating $\{x\}$ and $K$. By rotating $G$ we may w.l.o.g. assume that it lies parallel above the $d$-th basis axis with distance $c$. It follows that, using the convention that $\text{erfc}(0)=1$,

    $$\mu(K^c) \geq (2\pi)^{-1/2} \int\limits_{c} ^{\infty} \exp\left( -\frac{1}{2}  x^2 \right) \mathrm{d}x = \frac{1}{2} \text{erfc}(c/\sqrt{2}).$$
    Now, estimating the erfc term with a well-known lower bound (see e.g. \cite[inequality 7.1.13]{AbSt68}) 
    \begin{equation}
        \nonumber
        \begin{split}
            \frac{1}{2} \text{erfc}(c/\sqrt{2}) &> \frac{1}{\sqrt{\pi}} \frac{\exp\left( - \frac{1}{2} \log(\delta^{-1}) \right)}{\sqrt{\frac{\log(\delta^{-1})}{2}}+ \sqrt{\frac{\log(\delta^{-1})}{2} +2}} \\
            &\geq \frac{1}{\sqrt{\pi}} \frac{\delta^{1/2}}{2 \sqrt{\frac{\log(\delta^{-1})}{2} +2}} =:(*).
        \end{split}
    \end{equation}
    As
    $$\sqrt{\pi \left( 2 \log(\delta^{-1}) + 8 \right)} \leq \delta^{-1/2}$$
    by assumption,
    $$(*) \geq \frac{\delta^{1/2}}{\delta^{-1/2}} = \delta.$$
    Looking back, we have shown that 
    $$\mu(K^c) > \delta.$$
    The conclusion is that a set with $\mu(K) \geq 1 - \delta$ must contain all points in $B(\sqrt{\log(\delta^{-1})})$ as erfc$(\cdot)$ is non-increasing. This shows (\ref{ball_contained}).

    We denote 
    $$R:= \sqrt{ \log\left(\delta^{-1} \right)}, \: \: \:  \: \: \: C:= \frac{8 n^2}{n^2 -1} > 1$$
    and define
    $$d(x) := d(x, K).$$    
    With this notation the measure split is given by 
    \begin{equation}
    \label{31_claim}
        \begin{split}
            \nu_{bad}(\mathrm{d}x) \propto \mathrm{e}^{-\sigma ( d(x) )} \mu(\mathrm{d}x), \\
            \nu_{good}(\mathrm{d}x) \propto 1 - \mathrm{e}^{-\sigma ( d(x) )} \mu(\mathrm{d}x),
        \end{split}
    \end{equation}
    where $\sigma : [0,\infty] \rightarrow [0,R^2] $ is any $C^2 _b$ function satisfying
    \begin{enumerate}
        \item $\sigma$ is non-increasing,
        \item $\sigma(x) = R^2$ for $x \in [0,1]$,
        \item $\sigma(x) = 0$ for $x \geq 3 CR$,
        \item $|\sigma(x)'| \leq \frac{(x-1)}{C}$ for $x\geq 1$,
        \item $|\sigma(x)''| \leq \frac{1}{C}$ for $x\geq 0$.
    \end{enumerate}
    Such a function can be constructed explicitly (see \cite[Lemma 3.1]{Se22}). We now show that all claimed properties are fulfilled by the measures in (\ref{31_claim}).

    To show \ref{r1} simply estimate
    \begin{equation}
    \nonumber
    \begin{split}
        \delta' &= \expect \left[ \mathrm{e}^{-\sigma(d(x))}\right] \\
            &\leq \exp(-R^2) + 1 - \mu(K) \\
            &\leq \delta + \delta \\
            &= 2 \delta.
    \end{split}
\end{equation}

    For \ref{r2}, fix $\Tilde{C}>1$ and note that $d(x) \geq (\Tilde{C}-1)R$ whenever $x \notin \Tilde{C}K$, which shows that $d(x) \geq 3CR$ whenever $ x  \notin 4CK$. This implies that 
    $$\text{supp}(\nu_{\text{good}})\subseteq 4C K.$$

    By definition $K$ is symmetric and convex, so it holds that $\frac{\mathrm{d} \nu_{good}}{\mathrm{d} \mu}$ is (QC) and therefore \ref{r3} is fulfilled. 

    We finally show that 
    \begin{equation}
        \label{log_convave_goal}
        \frac{\mathrm{d}\nu_{bad}}{\mathrm{d} \mu^{\times n}} \propto \exp\left( - \sigma(d(x)) - \frac{n^2 - 1}{2 n^2} \lVert x \rVert^2 \right)
    \end{equation}
    is log-concave. The required symmetry such that the density is (QC) is evident. 

    To show convexity of 
    $$x \mapsto \sigma(d(x)) + \frac{n^2 - 1}{2 n^2} \lVert x \rVert^2$$
 take $x,z \in \rd$ and set for some $p \in (0,1)$
$$y := p x + (1-p)z .$$
Note that $$\sigma(d(x)) \geq \sigma ( \lVert x-P(y)\rVert ).$$
Of course, the same holds whenever $x$ is replaced by $z$. Using this inequality and rearranging 
$$\sigma(d(y)) + \frac{n^2 \! - \! 1}{2 n^2} \lVert y\rVert^2 \leq p \! \left( \! \sigma(d(x)) +  \frac{n^2 \! - \! 1}{2 n^2} \lVert x\rVert^2 \! \right) \! + (1-p)  \! \left( \! \sigma(d(z)) + \frac{n^2 \! - \! 1}{2 n^2} \lVert z\rVert^2 \! \right) \!,$$
the claim is proven whenever
\begin{equation}
\label{objective_equation_convexity}
    \begin{split}
        & p \sigma( \lVert  x - P(y)\rVert) + (1-p) \sigma( \lVert  z - P(y)\rVert) - \sigma ( d(y)) \\
        &\geq - \frac{n^2 - 1}{2 n^2} \left( p \lVert x\rVert^2 + (1-p) \lVert z \rVert^2 - \lVert y\rVert^2 \right) \\
        &= - p(1-p) \lVert x-z\rVert^2\frac{n^2 - 1}{2 n^2}
    \end{split}
\end{equation}
is shown. Denote the line connecting $x$ and $z$ by $\overline{xz}$ and let $o = p' x + (1-p')z$ be the projection of $P(y)$ onto $\overline{xz}$. Moreover, let 
$w_x := - \lVert  x- o\rVert$
and 
$w_z := \lVert  z - o \rVert.$
Clearly $|w_x - w_z| = \lVert x-z\rVert$.
By writing out the definitions,
$$ p w_x + (1-p) w_z = (p' -p)\lVert x-z\rVert.$$
It follows with
$w_y := \lVert  y - o \rVert = |p' - p | \:  \lVert x-z\rVert$
that 
$$w_y ^2  = \left( p w_x + (1-p) w_z \right)^2.$$
Define $a := \lVert  o - P(y)\rVert$ and assume w.l.o.g. that $a>0$. If this is not the case, do a limiting procedure in equation (\ref{convex_objective_2}) by using the continuity of $\sigma$. Then, 
$$d(x,P(y))^2 = a^2 + w_x ^2,$$
$$d(y,P(y))^2 = a^2 + \left( p w_x + (1-p) w_z \right)^2, $$
$$d(z,P(y))^2 = a^2 + w_z ^2.$$
Define 
$$f(w) := \sigma \left( \sqrt{a^2 + w^2} \right). $$
Using the chain rule,
$$f' (w) = \frac{w}{\sqrt{a^2 + w^2}} \sigma'\left(\sqrt{a^2 + w^2} \right),$$
$$f''(w) = \frac{w^2}{a^2 + w^2} \sigma'' \left(\sqrt{a^2 + w^2} \right) + \frac{a^2}{(a^2+w^2)^{3/2}} \sigma' \left(\sqrt{a^2 + w^2} \right).$$
Via property $5$ of $\sigma$,
$$\frac{w^2}{a^2 + w^2} \left|\sigma'' \left(\sqrt{a^2 + w^2} \right) \right| \leq \frac{1}{C}.$$
The fourth property allows to bound
$$\frac{a^2}{(a^2+w^2)^{3/2}} \left | \sigma' \left(\sqrt{a^2 + w^2} \right) \right| = \frac{a^2}{a^2 + w^2} \frac{ \left| \sigma' \left(\sqrt{a^2 + w^2} \right) \right|}{\sqrt{a^2 + w^2}} \leq \frac{1}{C}.$$
In total one obtains
\begin{equation}
    \label{sup_bound}
    \sup_{w \geq 0} |f''(w)| \leq \frac{2}{C}.
\end{equation}
The introduction of $f$ lets us rewrite the left-hand-side of equation (\ref{objective_equation_convexity}) to
\begin{equation}
    \label{convex_objective_2}
    \begin{split}
        & pf(w_x) + (1-p) f(w_z) - f(p w_x + (1-p)w_z).
    \end{split}
\end{equation}
A simple Taylor expansion gives (\ref{convex_objective_2}) $ \geq - 2 p(1-p) (w_x - w_z)^2 \sup\limits_{w \geq 0} | f''(w)|$. Applying estimate (\ref{sup_bound}) yields
$$pf(w_x) + (1-p) f(w_z) - f(p w_x + (1-p)w_z) \geq -\frac{ p(1-p) 4 \lVert x-z\rVert^2 }{C}.$$
By choice of $C$ inequality (\ref{objective_equation_convexity}) is fulfilled, which shows the claim.
    \end{proof}

\textbf{Acknowledgement:} we thank Volker Betz and Mark Sellke for helpful remarks. The author was partially supported by the  german research association (DFG) grant 535662048.

\end{document}